\documentclass{article}

\usepackage{proof}
\usepackage{latexsym}
\newtheorem{theorem}{Theorem}
\newtheorem{definition}[theorem]{Definition}
\newtheorem{proposition}[theorem]{Proposition}
\newtheorem{lemma}[theorem]{Lemma}
\newtheorem{cor}[theorem]{Corollary}
\newtheorem{mlemma}[theorem]{Main Lemma}
\newtheorem{claim}[theorem]{Claim}
\newcommand{\blem}{\begin{lemma}}
\newcommand{\elem}{\end{lemma}}
\newcommand{\bth}{\begin{theorem}}
\newcommand{\eth}{\end{theorem}}
\newcommand{\benu}{\begin{enumerate}}
\newcommand{\eenu}{\end{enumerate}}
\newcommand{\bdes}{\begin{description}}
\newcommand{\edes}{\end{description}}
\newcommand{\bdf}{\begin{definition}}
\newcommand{\edf}{\end{definition}}
\newcommand{\bcor}{\begin{cor}}
\newcommand{\ecor}{\end{cor}}
\newcommand{\bprp}{\begin{proposition}}
\newcommand{\eprp}{\end{proposition}}
\newcommand{\bmlem}{\begin{mlemma}}
\newcommand{\emlem}{\end{mlemma}}
\newcommand{\bclm}{\begin{claim}}
\newcommand{\eclm}{\end{claim}}
\newcommand{\bprf}{{\bf Proof}.\hspace{2mm}}

\newcommand{\eprf}{\hspace*{\fill} $\Box$}

\newcommand{\beqn}{\begin{equation}}
\newcommand{\eeqn}{\end{equation}}
\newcommand{\beqnarr}{\begin{eqnarray}}
\newcommand{\eeqnarr}{\end{eqnarray}}
\newcommand{\beqnarrs}{\begin{eqnarray*}}
\newcommand{\eeqnarrs}{\end{eqnarray*}}

\newcommand{\spand}{\,\&\,}

\newcommand{\alp}{\alpha}
\newcommand{\eps}{\epsilon}
\newcommand{\veps}{\varepsilon}

\newcommand{\Del}{\Delta}
\newcommand{\ome}{\omega}
\newcommand{\Ome}{\Omega}
\newcommand{\bet}{\beta}
\newcommand{\gam}{\gamma}
\newcommand{\Gam}{\Gamma}

\newcommand{\sig}{\sigma}
\newcommand{\Sig}{\Sigma}
\newcommand{\tht}{\theta}

\newcommand{\Lam}{\Lambda}
\newcommand{\vphi}{\varphi}

\newcommand{\fal}{\forall}
\newcommand{\exi}{\exists}

\newcommand{\Rarw }{\Rightarrow}

\newcommand{\lrarw}{\leftrightarrow}
\newcommand{\Lrarw}{\Leftrightarrow}

\newcommand{\cala}{{\cal A}}
\newcommand{\calb}{{\cal B}}

\newcommand{\incl}{\subseteq}

\begin{document}

\title{Exact bounds on epsilon processes
}


\author{Toshiyasu Arai
\\
Graduate School of Science,
Chiba University
\\
1-33, Yayoi-cho, Inage-ku,
Chiba, 263-8522, JAPAN
}

\date{}

\maketitle

\begin{abstract}
In this paper we show that the lengths of the approximating processes in epsilon substitution method 
are calculable by ordinal recursions in an optimal way.
\end{abstract}

{\it Epsilon substitution method\/} is a method proposed by D. Hilbert to prove the consistency of (formal) theories.
The idea behind the method is that one could replace consistently
transfinite/non-computable objects as a figure of speech by finitary/computable ones 
as far as transfinite ones are finitely presented as axioms of a theory.
In other words,
the replacement ({\it epsilon substitution\/}) depends on contexts, i.e., formal proofs 
in which axioms for the transfinite objects occur.
If this attempt would be successfully accomplished, then the (1-)consistency of the theory follows.

For example, for first order arithmetic PA, replace each existential formula $\exi x F[x]$
by $F[\epsilon x.F[x]]$, where the {\it epsilon term\/} $\eps x.F[x]$ intends to denote
the least number satisfying $F[x]$ if such a number exists.
Otherwise it denotes an arbitrary object, e.g., zero.
Then PA is interpretable in an extended 'propositional calculus'  having the {\it epsilon axioms\/}:
\beqn\label{eq:WFep}
(\eps) \: F[t] \to \eps x. F[x]\not> t \land F[\eps x.F[x]]
\eeqn
The problem is to find a solving substitution which assigns numerical values to epsilon terms and under which all the epsilon axioms occurring in a given proof are true.

Hilbert's Ansatz is, starting with the null substitution $S^{0}$ which assigns zero to whatever, 
to approximate a solution by correcting false values step by step,
 and thereby generate the process $S^{0},S^{1},\ldots$ ({\it H-process\/}).
The problem is to show that the process terminates.

In \cite{esubjh}, \cite{esubid}, \cite{esubMahlo}, \cite{esubidea} and \cite{esubpi2},
we formulated H-processes for theories of jump hierarchies, 
for $ID_{1}(\Pi^{0}_{1}\lor\Sig_{1}^{0})$, for $[\Pi^{0}_{1},\Pi^{0}_{1}]$-FIX, for $\Pi^{0}_{1}$-FIX and for $\Pi^{0}_{2}$-FIX, resp., and proved that the processes terminate by transfinite induction up to the relevant proof-theretic ordinals.

In this paper we address a problem related to these termination proofs,
and show that the lengths of the processes are calculable by ordinal recursions in an optimal way.

Let T denote one of the following theories; first order arithmetic, the theories of
  the absolute jump hierarchy, theories $\Phi$-FIX for non-monotonic inductive definitions
 for the formula classes $\Phi=\Pi^{0}_{1}, [\Pi^{0}_{1},\Pi^{0}_{1}],\Pi^{0}_{2}$.
Let $|T|$ denote the proof-theoretic ordinal of T.

Given a finite sequence $Cr$ of critical formulas,
let $\{S^{n}\}$ denote the H-process for $Cr$.

\begin{theorem}\label{th:extbnd}
The length 
$H=\min\{n: S^{n} \mbox{ {\rm is a solution}}\}$ 
of the H-process up to reaching a solution for $Cr$,
is calculable by $|T|$-recursion.

Therefore so is the  solution $S^{H}$.
\eth

\section{First order arithmetic: Ackermann's proof}\label{sec:PA}

In this section we give the ordinal-theoretic heart of the epsilon substitution method.

\subsection{The H-process}

The language of first order arithmetic PA includes some symbols for computable functions, say $+$ for addition, $\cdot$ for multiplication and $\dot{-}$ for cut-off subtraction,
 and the relation symbol $<$. 
In its $\eps$-counterpart $\mbox{PA}\eps$, formulas and terms are defined simultaneously by stipulating that
\[
\mbox{if } F \mbox{ is a formula, then } \eps x. F \mbox{ is a term.}
\]
By {\it expression} we mean a term or a formula.

An $\eps$-{\it substitution\/} $S$ is a finite function assigning values $|\eps x. F|_{S}\in\ome$ of {\it canonical\/}(=closed and
minimal epsilon) terms $\eps x. F$. $dom(S)$ denotes its {\it domain\/}.

$\eps$-substitutions $S$ reduces an expression $e$ to its unique irreducible form $|e|_{S}$
by using default value $0$ for expressions not in $dom(S)$.

Let $Cr=\{Cr_{0},\ldots,Cr_{N}\}$ be a fixed finite sequence of closed epsilon axioms.
$S$ is {\it solving\/} if $S$ validates any critical formula in $Cr$. 
Otherwise $S$ is {\it nonsolving\/}.

The existence of a solving substitution for any finite sequence of critical formuls yields the 1-consistency of PA.

The {\it rank\/} $rk(e)<\ome$ of an expression $e$ measures nesting of bound variables in $e$.

\bdf
$rk(S):=\max(\{rk(e):e\in dom(S)\}\cup\{0\})$.
\edf

For a substitution $S$ and a natural number $r$, 
 $S_{<r}:=\{(e,v)\in S rk(e)<r\}$.

For a fixed sequence $Cr$, the {\it H-process\/} $S^{0}(=\emptyset),S^{1},\ldots$ of substitutions for $Cr$
 is defined using the ranks of $\eps$-terms.
The sequence $\{S^{n}\}$ is primitive (or even elementary) recursive.
We assume that  if $S^{n}$ is a solution for $Cr$,
then $S^{m}=S^{n}$ for any $m\geq n$.

 By an algorithm, we associate an epsilon axiom $Cr(S)$ to a nonsolving substitution $S$:
 \[
 Cr(S): \: F[t] \to \eps x. F[x]\not> t \land F[\eps x.F[x]],
 \]
 which is false under $S$.
 Then $e^{S}:\equiv \eps x. |F|_{S}$ and $v^{S}:=|t|_{S}$.
 
 If $S^{n}$ is nonsolving, then the next substitution is defined as follows.
\[
S^{n+1}:=S^{n}_{<rk(e^{S^{n}})}\cup\{(f,u)\in S^{n}: rk(f)=rk(e^{S^{n}}) \spand f\not\equiv e^{S^{n}}\}\cup\{(e^{S^{n}},v^{S^{n}})\}.
\]

\subsection{Termination proof}

In this subsection we recall a proof of the termination of the H-process.
The proof is based on the transfinite induction up to $\veps_{0}$.

Define the {\it Ackermann ordering\/}:
\beqn\label{eq:Ackorder}
x<_{A}y :\Lrarw [x\neq 0 \spand y=0] \, \lor\, [x,y\neq 0 \spand x<y]
\eeqn
Thus $0$ is the largest element in $<_{A}$.
$\|x\|_{A}$ denotes the order type of $x$ in the ordering $<_{A}$.

A relation $T\sqsubseteq_{A} S$ on $\eps$-substitutions is defined.

\bdf\label{df:sqsubset}
\beqnarrs
T\sqsubseteq_{A} S  & :\Lrarw & \fal (e,u)\in S\exi (e,v)\in T[v\leq_{A}u]  \\
& \Lrarw & |e|_{T}\leq_{A}|e|_{S} \mbox{ {\rm for any canonical} } e
\eeqnarrs
\edf

We associate an ordinal $ind(S)<\ome^{\ome}$ ({\it index\/} of $S$) relative to a fixed sequence $Cr$ of $\eps$-axioms.

$Cl_{\eps}(Cr)$ denotes the set of closed $\eps$-terms occurring in the set $Cr$. 
Let $N(Cr):=\# Cl_{\eps}(Cr)$(=the cardinality of the set $Cl_{\eps}(Cr)$). $N(Cr)$ is less than or equal to the total number of occurrences of the symbol $\eps$ in the set $Cr$.

\bdf\label{df:ordind}
\benu
\item\label{df:ordind.1}
{\rm For an} $e\in Cl_{\eps}(Cr)$ {\rm put}
\[\vphi(e;S):=\|v\|_{A} \mbox{ {\rm for }} v=|e|_{S}.\]
\item\label{df:ordind.2}
{\rm We arrange the set} $Cl_{\eps}(Cr)$ {\rm of cardinality} $N(Cr)$ {\rm as follows:} $Cl_{\eps}(Cr)=\{e_{i}:i<N(Cr)\}$ {\rm where} 
\[
e_{j} \mbox{ {\rm is a closed subexpression of }} e_{i} \Rarw j>i
\]

\item \label{df:ordind.3}
\[
ind(S)=\sum\{(\ome+1)^{i}\cdot\vphi(e_{i};S): i<N(Cr)\}.
\]
\eenu
\edf

Let $\mbox{IND}:=\mbox{IND}(Cr):=(\ome+1)^{N(Cr)}$.

Let $r_{n}=rk(S^{n})$, $e_{n}=e^{S^{n}}$, $v_{n}=v^{S^{n}}$ and $a_{n}=ind(S^{n})$ up to a solution. 
Otherwise let $r_{n}=e_{n}=v_{n}=a_{n}=0$.

The epsilon axiom $Cr(S)$ associated to nonsolving substitutions $S$ depends only on their indices $ind(S)$.

\blem\label{lem:PAB}(Cf. \cite{esubidea})\\
Let $S^{n}$ and $S^{m}$ be nonsolving substitutions such that $S^{m}\sqsubseteq_{A}S^{n}$. Then 
\benu
\item\label{lem:PAB.1} $a_{n}\geq a_{m}$. 
\item\label{lem:PAB.2} $S^{m+1}\sqsubseteq_{A}S^{n+1} \spand e_{n}=e_{m}\spand v_{n}=v_{m}$ and $r_{n+1}=r_{m+1}$ if $a_{n}= a_{m}$.
\eenu
\elem 

Each $S^{n}$ is shown to be {\it correct\/}, cf. \cite{esubidea}.
This yields the following fact for nonsolving $S^{n}$.
\beqn\label{eq:correctded}
(e_{n},v)\in S^{n} \Rarw 0\neq v_{n}<v
\eeqn

{\it Fix a positive integer\/} $\mbox{RANK}=\mbox{RANK}(Cr):=\max\{rk(Cr)+1,2\}$, where
$rk(Cr):=\max\{rk(Cr_{I}): I=0,\ldots,N\}$. 
Then for any $S$ appearing in the H-process, we have $rk(S)<\mbox{RANK}$.

Let
\[
\vec{S}^{m,k}=\{S^{n}\}_{m\leq n<  k}
.\]

\bdf\label{df:rkseries}
{\rm Let} $\vec{S}^{m,k}$ {\rm be a consecutive series in the H-process} $S^{0},\ldots$ {\rm Then}
\[
rk(\vec{S}^{m,k}):=\min(\{r_{i}: m<i< k\}\cup\{\mbox{{\rm RANK}}\})>0
.\]
\edf

\bdf\label{df:7.3}
{\rm A consecutive series} $\vec{S}^{m,k}$ {\rm in the H-process} $S^{0},\ldots$ {\rm is a} {\it section\/}  {\rm iff} 
$r_{m}<rk(\vec{S})$.
\edf

\bdf\label{df:73}
{\rm Let} $\vec{S}^{i}=\vec{S}^{m^{i},k^{i}}\, (i=0,1)$
{\rm be two consecutive series in the H-process} $S^{0},\ldots$
{\rm such that} $r_{m^{i}}\leq rk(\vec{S}^{i})$ {\rm for} $i=0,1$.

{\rm If} 
$S^{m^{1}}\sqsubseteq_{A}S^{m_{0}}$ {\rm and one of the following conditions is fulfilled, then we write} 
$\vec{S}^{1}\prec\vec{S}^{0}${\rm :}
 \benu
 \item {\rm There exists a} $p<\min\{\ell^{0},\ell^{1}\}\, (\ell^{i}:=k^{i}-m^{i})$ {\rm such that} $a_{m^{0}+p}>a_{m^{1}+p}$ {\rm and}
 $\fal i<p(a_{m^{0}+i}=a_{m^{1}+i})$.
 \label{lem:73.2a}
 \item $\ell^{1}<\ell^{0}$ {\rm and} 
 $\fal i<\ell^{0}(a_{m^{0}+i}=a_{m^{1}+i})$.
 \label{lem:73.2b}
 \eenu
\edf

The following Lemma \ref{lem:PAC} is seen readily from Lemma \ref{lem:PAB} and (\ref{eq:correctded}), cf. \cite{esubidea}.

\blem\label{lem:PAC}
 Let $\vec{S}^{i}=\vec{S}^{m^{i},k^{i}}\, (i=0,1)$ be two sections in the H-process $S^{0},\ldots$ such that 
 $k^{0}=m^{1}$ and $r_{m^{0}}\leq r_{m^{1}}<rk(\vec{S}^{0})$.
 Then
\benu
\item \label{lem:PAC.1}
$S^{m^{1}}\sqsubseteq_{A}S^{m^{0}}$.
\item \label{lem:PAC.2}
$\vec{S}^{1}\prec\vec{S}^{0}$.
\eenu
\elem

Lemma \ref{lem:PAC}.\ref{lem:PAC.2} means that each section $\vec{S}=\{S_{i}: i\leq k\}$ codes 
an ordinal $o(\vec{S})<\varepsilon_{0}$ in Cantor normal form with base $2$:
Let $rk(S_{0})\leq rk(\vec{S})=:r$.
Divide $\vec{S}$ into substrings which are sections as follows. Put $\{k_{0}<\cdots <k_{l}\}=\{i:i\leq k \spand rk(S_{i})=r\}\cup\{0\}$, and $\vec{S}=\vec{S}_{0}*\cdots*\vec{S}_{l}$ with $\vec{S}_{j}=(S_{k_{j}},\ldots,S_{k_{j+1}-1})$ for $0\leq j\leq l$ and $k_{l+1}=k+1$.

The series $\vec{S}_{0},\ldots,\vec{S}_{l}$ of substrings of $\vec{S}$ is called the {\it decomposition\/}\footnote{Note that the definition of the decomposition here differs from one in Definition \ref{df:decompose}.}
 of $\vec{S}$.

We have $\fal j<l[\vec{S}_{j+1}\prec\vec{S}_{j}]$.

For ordinals $a$ and $\alp\geq 2$ and $k<\ome$, let $\alp_{0}(a):=a$ and 
$\alp_{1+k}(a):=\alp^{\alp_{k}(a)}$.
Also set $\ome_{k}:=\ome_{k}(1)$.

For each series $\vec{S}=\vec{S}^{m,k}$ with $r_{m}\leq rk(\vec{S})$ 
and a natural number $\xi$ such that $0<\xi\leq r=rk(\vec{S})$,
associate an ordinal $o(\vec{S};\xi)<\ome_{\mbox{{\footnotesize RANK}}+2-\xi}$ so that the following Lemma \ref{lem:PAD}  holds,
cf. \cite{esubjh}  for a full definition and a proof.

\blem\label{lem:PAD} 
Let $\vec{S}^{i}=\vec{S}^{m^{i},k^{i}}\, (i=0,1)$ be two series in the H-process $S^{0},\ldots$ such that 
$\vec{S}^{1}\prec\vec{S}^{0}$ and $r_{m^{i}}<rk(\vec{S}^{i})$ for $i=0,1$.

Then
$o(\vec{S}^{0};\xi)>o(\vec{S}^{1};\xi)$ for any natural number $\xi\leq \min\{rk(\vec{S}^{0}),rk(\vec{S}^{1})\}$.
\elem

\bth\label{th:6.13}(Transfinite induction up to $\varepsilon_{0}$)\\
The H-process $S^{0},\ldots$ terminates.
\eth
\bprf
Suppose the H-process $S^{0},\ldots$ is infinite and put $r_{n}=rk(S^{n})$.
 
Inductively we define a sequence $\{n_{i} :i\in\ome\}$ of natural numbers as follows. 
First set $n_{0}=0$. Suppose $n_{i}$ has been defined. Then put $\bet_{i}=\min\{r_{n}:n>n_{i}\}$ and $n_{i+1}=\min\{n>n_{i}: r_{n}=\bet_{i}\}$. 

Then Lemma \ref{lem:PAD} yields an infinite decreasing sequence of ordinals, viz. 
\\
$\fal i [o(\vec{S}_{i+1};\xi)<o(\vec{S}_{i};\xi)<\ome_{\mbox{{\footnotesize RANK}}}<\veps_{0}]$ for  
$\vec{S}_{i}=\vec{S}^{n_{i}, n_{i+1}}$ and $\xi=\bet_{0}+1\geq 2$.
\eprf
 
Therefore the H-process $S^{0},\ldots$ for any given sequence $Cr$ of critical formulas terminates. It provides a closed and solving substitution, which in turn yields the 1-consistency $\mbox{RFN}_{\Sig^{0}_{1}}(\mbox{PA})$ of PA stating that any PA-provable  $\Sig^{0}_{1}$-sentence
is true. 

However the above proof is not entirely satisfactory.
Specifically the 1-consistency of PA is known to be equivalent, over a weak arithmetic,
 to the principle $\mbox{PRWO}_{\veps_{0}}$, which says that there is no infinite primitive recursive descending chain of ordinals$<\veps_{0}$,
or to be equivalent to the totality of  $\veps_{0}$-recursive functions.
The sequence $\{n_{i}\}_{i}$ and hence the sequence $\{o(\vec{S}_{i};\xi)\}_{i}$  of ordinals
 in the above proof are not seen to be recursive.
Therefore we need to show that the sequence $\{n_{i}\}_{i}$ is $\veps_{0}$-recursive
in showing the 1-consistency of PA.

\section{Exact bounds: finite ranks}

In this section we show that the length of the H-process up to reaching a solution is bounded by 
an ordinal recursive function.
From the bound one can easily read off the bound for the provably recursive functions in PA.

\subsection{Ordinal recursive functions}

Let us recall the definition and facts on ordinal recursive functions in W. W. Tait\cite{tait}.

Let $<_{\Lam}$ denote a primitive recursive well ordering of type $\Lam>0$.
Assume that $0$ is the least element in $<_{\Lam}$.

For each $\alp\leq\Lam$, $<_{\alp}$ denotes the initial segment of $<_{\Lam}$ of type $\alp$.
A number-theoretic function is said to be $\alp${\it -recursive\/} iff
it is generated from the schemata for primitive recursive functions plus the following schema
for introducing a function $f$ in terms of functions $g,h$ and $d$:
\[
f(\vec{y},x)=\left\{
\begin{array}{ll}
g(\vec{y},x) & \mbox{ if } d(\vec{y},x)\not<_{\alp}x  \\
h(\vec{y},x,f(\vec{y},d(\vec{y},x))) & \mbox{ if } d(\vec{y},x)<_{\alp}x 
\end{array}
\right.
\]
A function is $<\!\alp${\it -recursive\/} iff it is $\bet$-recursive for some $\bet<\alp$.

W. W. Tait\cite{tait}, p.163 shows that for each $\alp$ the class of $\alp$-recursive functions is closed under the
{\it external recursion\/} to introduce a function $f$ in terms of functions $g,h,d$ and $e$:
\[
f(\vec{y},x)=\left\{
\begin{array}{ll}
g(\vec{y},x) & \mbox{ if } e(\vec{y},d(\vec{y},x))\not<_{\alp}e(\vec{y},x)  \\
h(\vec{y},x,f(\vec{y},d(\vec{y},x))) & \mbox{ if } e(\vec{y},d(\vec{y},x))<_{\alp}e(\vec{y},x)
\end{array}
\right.
\]

\subsection{$p$-series}

In this subsection we define a series $\vec{S}^{m,k}$ to be a $p${\it -series\/}.
$p$-series is introduced for counting the number of ranks $r_{n}$ in the H-process.koko

Given the finite sequence $Cr=\{Cr_{I}: I\leq N\}$ of critical formulas in $\mbox{PA}\eps$, 
let $\{S^{n}\}$ denote the H-process for $Cr$.
Recall that the sequence is infinite in the sense that if $S^{n}$ is a solution for $Cr$,
then $S^{m}=S^{n}$ for any $m\geq n$.

Recall that $\ome^{\ome}>\mbox{IND}=\mbox{IND}(Cr):=(\ome+1)^{N(Cr)}>a_{n}$ and 
$\ome>\mbox{RANK}=\mbox{RANK}(Cr)>r_{n}$
for any $n$.

For $m<k$ let
\beqnarrs
{\sf nd}(\vec{S}^{m,k}) & := & \{n\in [m,k): r_{n}\leq rk(\vec{S}^{n,k})\}
\\
& (=& \{n\in [m,k): r_{n}=\min(r_{i}: i\in [n,k))\}).
\eeqnarrs

\bdf\label{df:decompose}
{\rm Let}  $\vec{S}=\vec{S}^{m,k}$ {\rm (with} $m<k${\rm ) such that} 
$r_{m}\leq rk(\vec{S})${\rm (i.e.,} $m\in{\sf nd}(\vec{S})${\rm ),
and let}
$\{k_{0},\ldots,k_{l}\}_{<}={\sf nd}(\vec{S})${\rm . Then}
$(\vec{S}_{0},\ldots,\vec{S}_{l})$ {\rm with} $\vec{S}_{j}:=\vec{S}^{k_{j},k_{j+1}}$ {\rm and} $k_{l+1}:=k$
{\rm is called} {\it the decomposition of\/} $\vec{S}$ {\it into substrings\/}.
{\rm Each substring} $\vec{S}_{j}\, (0\leq j\leq l)$ {\rm is called a} {\it component\/} {\rm in the decomposition of} $\vec{S}$
\edf

Note that $k_{0}=m$, $k_{l}=k-1$, and
$rk(S^{k_{j}})\leq rk(S^{k_{j+1}})<rk(\vec{S}_{j})$ for $j<l$.
Also note that each component $\vec{S}_{j}$ is a section.

\blem\label{lem:concatenation}
  Let $\vec{S}^{i}=\vec{S}^{m^{i},k^{i}}\, (i=0,1)$ with $k^{0}=m^{1}$
 such that $r_{m^{i}}\leq rk(\vec{S}^{i})\, (i=0,1)$ and $r_{m^{0}}\leq r_{m^{1}}$.
 Then for $\vec{S}:=\vec{S}^{0}*\vec{S}^{1}$ we have
 \[
 {\sf nd}(\vec{S})=\{n\in{\sf nd}(\vec{S}^{0}): n\leq I\}\cup{\sf nd}(S^{1})
 ,\]
 where $I:=\max\{n\in [m^{0},k^{0}): r_{n}\leq r_{m^{1}}\}$.
\elem
\bprf
We see $I=k^{0}_{J}$ for a $J\leq l^{0}$ from the facts that both $\{S^{n}: I\leq n<m^{1}\}$
and each $\vec{S}^{0}_{j}\, (j\leq l^{0})$
 are sections and $k^{0}_{l^{0}}=m^{1}-1$.
Therefore $m^{1}=k^{1}_{0}=k_{J+1}$ and the lemma is shown.
\eprf

\bdf\label{df:pseries}
{\rm Let} $\vec{S}=\vec{S}^{m,k}$ {\rm with} $m<k$ {\rm
such that}  $r_{m}\leq rk(\vec{S})${\rm . 
Define inductively the series} $\vec{S}$  {\rm to be a}  $p${\it -series\/} {\rm and a} $p${\it -section\/} {\rm as follows:}
\benu
\item
$\vec{S}$ {\rm is a} $0${\it -series\/} {\rm iff} $k=m+1${\rm , i.e., a singleton.}
\item {\rm A}  $p${\rm -series is  a}  $p${\it -section\/} {\rm iff it is a section.}
\item
{\rm Let}  $\vec{S}=\vec{S}_{0}*\cdots*\vec{S}_{l}$
{\rm be the decomposition of} $\vec{S}$ {\rm into substrings. Then}
$\vec{S}$  {\rm is a}  $(p+1)${\it -series\/} {\rm iff
each substring} $\vec{S}_{j}$ {\rm is  a}  $p${\rm -section, or equivalently a} $p${\rm -series.}
\eenu
\edf

 \blem\label{lem:pseries}
 \benu
 \item\label{lem:pseries.1}
 Each $p$-series is a $(p+1)$-series.

 \item\label{lem:pseries.3}
 Let $\vec{S}^{i}=\vec{S}^{m_{i},k_{i}}\, (i=0,1)$ be two $p$-series overlapped, i.e., $[m_{0},k_{0})\cap[m_{1},k_{1})\neq\emptyset$.
 Then the union  $\vec{S}=\vec{S}^{\min\{m_{0},m_{1}\},\max\{k_{0},k_{1}\}}$
 is a $p$-series.
 
  \item\label{lem:pseries.4}
  Let us call a $p$-series {\rm proper} if $p=0$, or $p>0$ and it is not a $(p-1)$-series.
 \benu
 \item\label{lem:pseries.41}
 If $\vec{S}$ is a proper $p$-section, then
  \[
  \#\{rk(S): S\in\vec{S}\}\geq p+1.
  \]
  \item\label{lem:pseries.42}
   If $\vec{S}$ is a proper $p$-series, then
  \[
 \#\{rk(S): S\in\vec{S}\}\geq p.
  \]
   \item\label{lem:pseries.43}
  If a proper $p$-series $\vec{S}$ begins with $S^{0}=\emptyset$, then
    \[
  \#\{rk(S): S\in\vec{S}\}\geq p+1.
  \]
  Therefore there is no proper {\rm RANK}-series beginning with $S^{0}$. 
  \eenu
 \eenu
 \elem
 \bprf
 By induction on $p$.
 \\
\ref{lem:pseries}.\ref{lem:pseries.1}.
 A  $0$-series $\{S^{n}\}$ is a $1$-series.
 \\
 \ref{lem:pseries}.\ref{lem:pseries.3}.
Assume $p>0$ and one is not a substring of the other, i.e., $[m_{i},k_{i})\not\incl[m_{1-i},k_{1-i})$.
Then without loss of generality we may assume $m_{0}<m_{1}<k_{0}<k_{1}$.
Decompose the $p$-series $\vec{S}^{i}$ to the sequence of $(p-1)$-series
$\vec{S}^{i}_{j}=(S^{k^{i}_{j}},\ldots,S^{k^{i}_{j+1}-1})\, (j\leq l_{i})$.
It suffices to show:
$m_{1}\leq k_{j}^{i}\leq k_{0} \Rarw 
 \exi j^{\prime}(k^{i}_{j}=k^{1-i}_{j^{\prime}})$.
 
This is seen from the condition that each decomposition $\{\vec{S}^{i}_{j}: j\leq l_{i}\}$ is a 
sequence of sections  with nondecreasing ranks of the first terms.
 \\
\ref{lem:pseries}.\ref{lem:pseries.4}.
Let  $\vec{S}=\vec{S}_{0}*\cdots*\vec{S}_{l}$ be a proper $(p+1)$-series 
with ${\sf nd}(\vec{S})=\{k_{0},\ldots,k_{l}\}_{<}$. 
Then $l>0$ and one of $p$-sections $\vec{S}_{j}$ is proper.
Lemma \ref{lem:pseries}.\ref{lem:pseries.41} yields
$\#\{rk(S): S\in\vec{S}_{j}\}\geq p+1$, and hence Lemma \ref{lem:pseries}.\ref{lem:pseries.42} 
follows.

If $\vec{S}_{0}$ is proper, then $r_{k_{0}}<r_{k_{1}}<rk(\vec{S}_{0})$
since $\vec{S}_{0}$ is a section. Hence $\#(\{rk(S): S\in\vec{S}_{0}\}\cup\{r_{k_{1}}\})\geq p+2$.
Next assume $j>0$. Then $r_{k_{0}}<r_{k_{j}}<rk(\vec{S}_{j})$, and
$\#(\{rk(S): S\in\vec{S}_{j}\}\cup\{r_{k_{0}}\})\geq p+2$.
This shows Lemma \ref{lem:pseries}.\ref{lem:pseries.41}.

Lemma \ref{lem:pseries}.\ref{lem:pseries.43}  is seen from the fact $r_{n}>0$ for $n>0$.
Namely any proper $p$-series $\vec{S}$ beginning with $S^{0}=\emptyset$ is a section.
 \eprf

\blem\label{lem:opseries.1}
  Let $\vec{S}^{i}=\vec{S}^{m^{i}, k^{i}}\, (i=0,1)$ be two consecutive series, $k^{0}=m^{1}$
 such that $r_{m^{i}}\leq rk(\vec{S}^{i})\, (i=0,1)$ and $r_{m^{0}}\leq r_{m^{1}}$.

The concatenated series $\vec{S}=\vec{S}^{m_{0}, k_{1}}$ is a $(p+1)$-series
if $\vec{S}^{0}$ is a $p$-series and $\vec{S}^{1}$ is a $(p+1)$-series.
\elem
\bprf
This is seen from Lemmas \ref{lem:concatenation} and \ref{lem:pseries}.\ref{lem:pseries.1}.
\eprf

Let $<_{\veps_{0}}$ denote a standard well ordering of type $\veps_{0}$ with the least element $0$.

\blem\label{lem:opseries}
Let $\vec{S}^{i}=\vec{S}^{m^{i}, k^{i}}\, (i=0,1)$ be two $p$-series such that $k^{0}=m^{1}$,
$S^{k^{1}-1}$ is nonsolving and $r_{m^{0}}\leq r_{m^{1}}$.
For $o(\vec{S}^{i}):=o(\vec{S}^{i};0)$ we have 
$o(\vec{S}^{1})<_{\veps_{0}}o(\vec{S}^{0})$.
\elem
\bprf
By Lemma \ref{lem:PAD} it suffices to show $\vec{S}^{1}\prec\vec{S}^{0}$.
As in Lemma \ref{lem:PAC} this is seen as follows.
Let $\ell^{i}:=k^{i}-m^{i}$.

Since the relation $\sqsubseteq_{A}$ is transitive, we have $S^{m^{1}}\sqsubseteq_{A} S^{m^{0}}$
by Lemma \ref{lem:PAC}.\ref{lem:PAC.1}.
Using Lemma \ref{lem:PAB}, it suffices to show that the following case never happen: 
$\ell^{0}\leq\ell^{1}$
and 
$\fal i<\ell^{0}[a_{m^{0}+i}=a_{m^{1}+i}]$.

 
If this happens, then we would have 
$a_{m^{1}-1}=a_{m^{0}+\ell^{0}-1}=a_{m^{1}+\ell^{0}-1}$, and hence
$S^{m^{1}}\ni(e_{m^{1}-1},v_{m^{1}-1})=(e_{m^{1}+\ell^{0}-1},v_{m^{1}+\ell^{0}-1})$ by Lemma \ref{lem:PAB}.\ref{lem:PAB.2}.
On the other hand we have
$S^{m^{1}+\ell^{0}-1}\sqsubseteq_{A}S^{m^{1}}$ by Lemma \ref{lem:PAC}.\ref{lem:PAC.1}.

For a $v\leq_{A}v_{m^{1}+\ell^{0}-1}$ we would have $(e_{m^{1}+\ell^{0}-1},v)\in S^{m^{1}+\ell^{0}-1}$.
By (\ref{eq:correctded}) we have $0\neq v_{m^{1}+\ell^{0}-1}<v$.
A contradiction.
\eprf

$k=M(p,n)$ defined below will denote the number such that $\vec{S}^{n,k}$ is the longest
$p$-series starting  with nonsolving $S^{n}$.

\bdf\label{df:Lp}

$M(0,n):=n+1$.

 \bdes
 \item[Case 0] $S^{n}$ {\rm is solving:} $M(p+1,n):=n$.
 
 \item[Case 1] $S^{n}$ {\rm is nonsolving. Let}
 \[
 e_{p}(n):=o(\vec{S}^{n,M(p,n)})
 .\]
 {\rm Then define}
 \[
 M(p+1,n):=\left\{
 \begin{array}{ll}
 M(p,n) & \mbox{{\rm if }} e_{p}(M(p,n))\not<_{\veps_{0}}e_{p}(n)
 \\
  M(p,n) & \mbox{{\rm if }} r_{M(p,n)}<r_{n} \spand e_{p}(M(p,n))<_{\veps_{0}}e_{p}(n)
  \\
M(p+1,M(p,n)) & \mbox{{\rm if }} r_{M(p,n)}\geq r_{n} \spand e_{p}(M(p,n))<_{\veps_{0}}e_{p}(n)
 \end{array}
 \right.
 \]
 \edes
 \edf

Actually the function $M(p,n)$ depends also on the given sequence  $Cr$ of epsilon axioms.
We  write $M(p,n;Cr)$ for $M(p,n)$ when the parameter  $Cr$ should be mentioned.

  
 A consecutive series $\vec{S}^{n,k}$ is a {\it normal\/} $p${\it -series\/} iff it is a $p$-series and $S^{k-1}$ is nonsolving if $k>n$.
 
 \blem\label{lem:Lp}
 \benu

   \item\label{lem:Lp.1}
   If $S^{n}$ is nonsolving, then $\vec{S}^{n,M(p,n)}$ is a normal $p$-series.

 \item\label{lem:Lp.2}
 If $\vec{S}^{n,k}$ is a normal $p$-series, then $k\leq M(p,n)$.

  \item\label{lem:Lp.3}
  $S^{H}$ is a solution for $Cr$, where
  $H=H(Cr):=M(\mbox{{\rm RANK}}-1,0;Cr)$.

 \eenu
 \elem
 \bprf
\\
 \ref{lem:Lp}.\ref{lem:Lp.1}. Main induction on $p$.
 The case when $p=0$ is trivial.
 
 The case $p+1$ is proved by side induction on $e_{p}(n)$.
 Assume that $S^{n}$ is nonsolving.
 \benu
 \item
 $M(p+1,n)=M(p,n)$:
 Then by MIH(=Main Induction Hypothesis) $\vec{S}^{n,M(p,s)}$ is a normal $p$-series.
 $\vec{S}^{n,M(p,s)}$ is also a normal $(p+1)$-series by Lemma \ref{lem:pseries}.\ref{lem:pseries.1}.
 \item
$M(p+1,n)\neq M(p,n)$:
 Then with $k=M(p,n)$ we have $r_{k}\geq r_{n}$ and $M(p+1,n)=M(p+1,k)$.
 By MIH $\vec{S}^{n,k}$ is a normal $p$-series.
 Since $M(p+1,k)\neq k$, $S^{k}$ is nonsolving, and hence by MIH, $\vec{S}^{k,M(p,k)}$ is a normal $p$-series.
 Lemma \ref{lem:opseries} yields $e_{p}(n)=o(\vec{S}^{n,k})>o(\vec{S}^{k,M(p,k)})=e_{p}(k)$.
 Therefore $\vec{S}^{k,M(p+1,k)}$ is a normal $(p+1)$-series by SIH(=Side Induction Hypothesis).
 Together with $r_{n}\leq r_{k}$ it follows from  Lemma \ref{lem:opseries.1} that
 $\vec{S}^{n,M(p+1,n)}$ is a normal $(p+1)$-series.
\eenu
\ref{lem:Lp}.\ref{lem:Lp.2}. Main induction on $p$.
   The case when $p=0$ is trivial.
   
The case $p+1$. First we show the following:
\beqn\label{lem:Lp.15}
  n\leq n^{\prime}< M(p,n) \Rarw M(p,n^{\prime})\leq M(p,n)
\eeqn
 Assume $n\leq n^{\prime}< M(p,n)=:k$ and $n^{\prime}<M(p,n^{\prime}):=k^{\prime}$.
 Then by Lemma \ref{lem:Lp}.\ref{lem:Lp.1} $\vec{S}^{n,k}$ and $\vec{S}^{n^{\prime},k^{\prime}}$ are two normal $p$-series overlapped.
 By Lemma  \ref{lem:pseries}.\ref{lem:pseries.3} the union $\vec{S}^{n,\max\{k,k^{\prime}\}}$ is a normal $p$-series too.
 By MIH it follows that $k^{\prime}\leq M(p,n)$.
 This shows (\ref{lem:Lp.15}).
 
 Now by side induction on $k-n$ we prove:
 \[
 \mbox{If } \vec{S}^{n,k} \mbox{ is a normal } (p+1)\mbox{-series, then } k\leq M(p+1,n)
.\]
Assume that  $\vec{S}^{n,k}$ is a normal $(p+1)$-series, and ${\sf nd}(\vec{S}^{n,k})=\{k_{0},\ldots,k_{l}\}_{<}$, $l>0$.
Then by MIH we have $k_{1}\leq M(p,n)$.
Let $j\leq l$ denote maximal such that $k_{j}\leq M(p,n)$.
\benu
\item
$k_{j}=M(p,n)$: 
Since $\vec{S}^{k_{j},k}$ is a normal $(p+1)$-series, we have by SIH that $k\leq M(p+1,k_{j})$.
On the other hand we have $M(p+1,n)=M(p+1,M(p,n))$ by Definition \ref{df:Lp}, $r_{k_{j}}\geq r_{k_{0}}=r_{n}$ and
$e_{p}(M(p,n))<_{\veps_{0}}e_{p}(n)$, Lemma \ref{lem:opseries}.
Hence $k\leq M(p+1,k_{j})=M(p+1,n)$.

\item $k_{j}<M(p,n)$: 
 \benu
 \item $j=l$:
 Then $k_{l}<M(p,n)$, and hence $k=k_{l}+1\leq M(p,n)\leq M(p+1,n)$.
 
 \item $j<l$:
 Since $\vec{S}^{k_{j},k_{j+1}}$ is a normal $p$-series, we have $k_{j+1}\leq M(p,k_{j})$ by MIH.
 On the other hand we have $M(p,k_{j})\leq M(p,n)$ by (\ref{lem:Lp.15}).
 Thus $k_{j+1}\leq M(p,n)$, and this is not the case.
 \eenu
\eenu
 \ref{lem:Lp}.\ref{lem:Lp.3}.
 Let  $H=H(Cr):=M(p,0;Cr)$ for $p:=\mbox{{\rm RANK}}-1$.
 If $S^{0}$ is solving, then $0=H$.
 Suppose that $S^{0}$ is nonsolving.
 By Lemma \ref{lem:Lp}.\ref{lem:Lp.1} $\vec{S}^{0,H}$ is a $p$-series.
 From Lemma \ref{lem:opseries.1} and $r_{0}=0\leq r_{H}$ we see that 
 $\vec{S}^{0,H+1}$ is a $(p+1)$-series.
But this means that  $\vec{S}^{0,H+1}$ is a $p$-series by Lemma \ref{lem:pseries}.\ref{lem:pseries.43}.
Therefore we see from Lemma \ref{lem:Lp}.\ref{lem:Lp.2} that $\vec{S}^{0,H+1}$ is not normal, i.e., $S^{H}$ is solving.
 \eprf

  \blem \label{lem:nested}
 The function  $(p,n,Cr)\mapsto M(p,n;Cr)$ is $\veps_{0}$-recursive.
 \elem
 \bprf
It suffices to see that  $M(p,n;Cr)$ is defined by nested recursion on the ordinal $\veps_{0}$.
Then it is $\veps_{0}$-recursive by a result in W. W. Tait\cite{nested} and $\ome^{\veps_{0}}=\veps_{0}$.

Suppressed the parameter $Cr$, let us define a function $M^{\prime}(p,n,y)$ as follows:

$M^{\prime}(0,n,y):=n+1$.

 \bdes
 \item[Case 0] $S^{n}$ {\rm is solving:} $M^{\prime}(p+1,n,y):=n$.
 
 \item[Case 1] $S^{n}$ {\rm is nonsolving.}
  \benu
  \item $o(\vec{S}^{M^{\prime}(p,n,y),M^{\prime}(p,M^{\prime}(p,n,y),y)})\not<_{\veps_{0}}y$:
 \[
 M^{\prime}(p+1,n,y):=M^{\prime}(p,n,y)
 .\]
 \item $ r_{M^{\prime}(p,n,y)}<r_{n} \spand 
  o(\vec{S}^{M^{\prime}(p,n,y),M^{\prime}(p,M^{\prime}(p,n,y),y)})<_{\veps_{0}}y$:
  \[
  M^{\prime}(p+1,n,y):= M^{\prime}(p,n,y)
  .\]
  \item
  $r_{M^{\prime}(p,n,y)}\geq r_{n} \spand 
 o(\vec{S}^{M^{\prime}(p,n,y),M^{\prime}(p,M^{\prime}(p,n,y),y)})<_{\veps_{0}}y$:
   \[
  M^{\prime}(p+1,n,y):= M^{\prime}(p+1,M^{\prime}(p,n,y), o(\vec{S}^{M^{\prime}(p,n,y),M^{\prime}(p,M^{\prime}(p,n,y),y)}))
  .\]
  \eenu
 \edes
Then $M^{\prime}(p,n,y)$  is seen to be defined by nested recursion on the lexicographic ordering $\prec$ on pairs
$(p,y)$:
$(p,y)\prec(q,z)$ iff $p<q$ or $p=q \spand y<_{\veps_{0}}z$.

Then $M(p,n):=M^{\prime}(p,n,\ome_{\mbox{{\footnotesize RANK}}+2})$ enjoys the defining clauses in Definition \ref{df:Lp}.
 \eprf

 \section{Exact bounds: infinite ranks}
 
 In this section let us compute the length of the H-process for theories of jump hierarchies, which is slightly modified from \cite{esubjh}.

The normal function $\tht_{\alp}: \bet\mapsto \tht\alp\bet$ is the $\alp^{\mbox{{\small th}}}$ iterate of the function $\tht 1\bet=\ome^{\bet}$.
 Fix an ordinal $\Lam<\Gam_{0}$, the least strongly critical number,
and let $<_{\Gam_{0}}$ denote a standard primitive recursive well ordering of type 
$\Gam_{0}$ with the least element $0$.
In what follows the subscript in $<_{\Gam_{0}}$ is omitted.

Let $\cala(x,\alp,z,X)$ be a fixed quantifier free formula in the language of first order arithmetic.
Let $(H)_{\Lam}$ denote the theory of
  the absolute jump hierarchy $\{H_{\alp}\}_{\alp\leq\Lam}$ 
  generated by the formula $\cala$ and up to $\alp\leq\Lam$:
\[
\alp\leq_{\Gam_{0}}\Lam \to \{y\in H_{\alp} \lrarw \exi x\cala(x,\alp,y,H_{<\alp})\}
\]
where $H_{<\alp}=\sum_{\bet<\alp}H_{\bet}$, i.e., $H_{<\alp}$ denotes the binary abstract $\{(\bet,z): (\bet,z)\in H_{<\alp}\}$.

 The critical formulas in its $\eps$-counterpart are the $\eps$-axiom (\ref{eq:WFep}),
 \[
 \alp\leq\Lam \to 
 \{t\in H_{ \alp}\lrarw \cala(\eps x\cala(x, \alp,t,H_{< \alp}), \alp,t,H_{< \alp})\}
 \]
 and
 \[
 \alp\leq\Lam \to \{(\bet,t)\in H_{< \alp}\lrarw \bet<  \alp \land t\in H_{\bet}\}
 \]

Now the rank of an expression is defined such that $rk(e)<3\Lam+\ome$.
 For a given finite sequence $Cr$ of critical formulas, 
 let $\mbox{RANK=RANK}(Cr)=3\Lam+n>rk(Cr):=\max\{rk(Cr_{I}): I=0,\ldots,N\}$ for an $n<\ome$.

An $\eps$-substitution is a finite function assigning numerical 
values $|\eps x. F|_{S}\in\ome$ to canonical terms $\eps x. F$, 
and boolean values $|e|_{S}=\top$ for expressions $e$ in one of the shapes
$n\in H_{\alp}$ or $(\bet,n)\in H_{<\alp}$ such that $\bet<\alp\leq\Lam$.

Define the Ackermann ordering $<_{A}$ and $\|x\|_{A}$ as in (\ref{eq:Ackorder}),
where for boolean values, $\bot<_{A}\top$ and $\|\bot\|_{A}=0, \|\top\|_{A}=1$.

Define the index $ind(S)<\mbox{IND}=\mbox{IND}(Cr)=(\ome+1)^{N(Cr)}<\ome^{\ome}$ of $S$ relative to a fixed sequence $Cr$ of critical formulas as in Definition \ref{df:ordind}.

 \subsection{Bounds on $p$}
 
 Since $\mbox{RANK}\geq\ome$, Lemma \ref{lem:pseries}.\ref{lem:pseries.43}
 is useless here.
  We need to give a bound on $p$ such that $M(p,n;Cr)<M(p+1,n;Cr)$.
  
 Let $\ell_{p}=M(p+1,0;Cr)$. $\{\ell_{p}\}$ is an increasing sequence $\ell_{p}\leq\ell_{p+1}$,
 and once $\ell_{p}=\ell_{p+1}$, then $\ell_{p}=\ell_{q}$ for any $q\geq p$.
 
 Now let $p(Cr)$ denote the least number $p$ such that $\ell_{p}=\ell_{p+1}$ 
 if such a $p$ exists.
 We show that the number $p(Cr)$ is defined.
 Then $S^{H}$ is a solution of $Cr$ for $H=M(p(Cr),0;Cr)$.
 
For each $p$,
let $\vec{S}_{p+1}=\{S^{n}\}_{n<\ell_{p}}$ be 
the $(p+1)$-section according to Lemma  \ref{lem:Lp}.\ref{lem:Lp.1}, 
and $\vec{S}_{p+1}=\vec{S}_{p}^{1}*\cdots*\vec{S}_{p}^{l_{p}}$ 
its decomposition into $p$-sections $\vec{S}_{p}^{j}$.

Let $o(\vec{S})=o(\vec{S};0)<\tht(\mbox{RANK})(\mbox{IND})<\tht\Lam\veps_{0}$ 
denote the ordinal associated to sections $\vec{S}$ as in \cite{esubjh}.

Then let $\alp_{p}^{j}=o(\vec{S}_{p}^{j})$ for $0<j< l_{p}$,
and $\alp_{p}^{0}:=\tht(\mbox{RANK})(\mbox{IND})$.

By Lemma \ref{lem:opseries}  we have
\beqn\label{eq:alpdec}
\alp_{p}^{j}>\alp_{p}^{j+1}
\eeqn

Moreover let $\gam_{p}^{j}=r_{k_{p}^{j+1}}\, (0\leq j< l_{p})$,
where $S^{k_{p}^{j}}$ is the first term in the substring $\vec{S}_{p}^{j+1}$ for $j<l_{p}$, and
$k_{p}^{l_{p}}=\ell_{p}$.

Finally let for $\Del:=\mbox{RANK}>\gam_{p}^{j}$
\[
\tht\Lam\veps_{0}>\bet_{p}  :=  \sum_{j<l_{p}}\Del^{\alp_{p}^{j}}\gam_{p}^{j} 
= \Del^{\alp_{p}^{0}}\gam_{p}^{0}+\cdots+\Del^{\alp_{p}^{l_{p}-1}}\gam_{p}^{l_{p}-1}
\]

From (\ref{eq:alpdec}) we see that $\bet_{p}$ is in Cantor normal form.

\blem\label{lem:pbnd}
$\bet_{p}>\bet_{p+1}$ 
if $\ell_{p}<\ell_{p+1}$.
\elem
\bprf
Assume $\ell_{p}<\ell_{p+1}$.
This means that $S^{\ell_{p}}$ is nonsolving.

By Lemma  \ref{lem:Lp}.\ref{lem:Lp.2}  we have 
$\gam_{p}^{l_{p}-2}>\gam_{p}^{l_{p}-1}$.

Let
\[
J=\min\{j\leq l_{p}-2 : \gam_{p}^{j}>\gam_{p}^{l_{p}-1}\}
.\]

Then 
\[
k_{p}^{J}=\max\{i<\ell_{p}: r_{i}\leq \gam_{p}^{l_{p}-1}\}
\]
and by Lemma \ref{lem:concatenation} we have $\vec{S}_{p+1}^{j}=\vec{S}_{p}^{j}$ for $0<j\leq J$.
Hence $\alp_{p+1}^{j}=\alp_{p}^{j}$ for $0\leq j\leq J$, 
and $\gam_{p+1}^{j}=\gam_{p}^{j}$ for $0\leq j<J$.

Consider the next substitution $S^{\ell_{p}}$ to the last one in $\vec{S}_{p}^{l_{p}}$ or equivalently the last one in
$\vec{S}_{p+1}$.
Then $S^{\ell_{p}}$ is  the next substitution to the last one in $\vec{S}_{p+1}^{J+1}$
by Lemma \ref{lem:concatenation}.

Hence $\gam_{p}^{J}>\gam_{p}^{l_{p}-1}=\gam_{p+1}^{J}$.
Therefore $\bet_{p}>\bet_{p+1}$ as desired.
\eprf

Now $p(Cr)=F(0,Cr)$ for the function
\[
F(p,Cr)=\left\{
\begin{array}{ll}
F(p+1,Cr) & \mbox{ if }   \ell_{p+1}>\ell_{p}\\
\min\{q\leq p : \ell_{q+1}=\ell_{q}\} & \mbox{ otherwise}
\end{array}
\right.
\]

$F(p,Cr)$ is defined by a $\tht\Lam\veps_{0}$-external recursion
by Lemma \ref{lem:pbnd}.

 \section{Exact bounds: impredicative cases}
 
 Let T denote one of the theories $\Phi$-FIX for non-monotonic inductive definitions
 for the formula classes $\Phi=\Pi^{0}_{1}, [\Pi^{0}_{1},\Pi^{0}_{1}],\Pi^{0}_{2}$
(cf.  \cite{esubid}, \cite{esubMahlo}, \cite{esubidea} and \cite{esubpi2}).
 
 $O(\mbox{T})$ denotes the system of ordinal diagrams for T.
 T is a two sorted theory: one sort for natural numbers, and the other sort for ordinals.
 The well ordering $<$ on ordinals is understood to 
 be the ordering in the notation system $O(\mbox{T})$, 
 when the values of expressions are calculated.
 Its largest value is denoted $\pi\in O(\mbox{T})$, which is intended to be a closure
 ordinal of non-monotonic inductive definitions by the operators in $\Phi$.
 $\Ome\leq\pi$ is the first non-recursive ordinal $\ome_{1}^{CK}$.
 $d_{\Ome}\veps_{\pi+1}$ denotes the proof-theoretic ordinal of T,
 and the length of the H-process should be calculated by $d_{\Ome}\veps_{\pi+1}$-recursion.
 
The rank of an expression is defined such that $rk(e)<\pi+\ome$.
 For a given finite sequence $Cr$ of critical formulas, 
 let $\mbox{RANK=RANK}(Cr)=\pi+n(Cr)$ for an $n(Cr)<\ome$
so that
$\max\{rk(Cr_{I}): I=0,\ldots,N\}<\mbox{RANK}$. 
Then $\pi\neq rk(S)<\mbox{RANK}$ for any $S$ appearing in the H-process for $Cr$.

Define the index $ind(S)<\mbox{IND}=\mbox{IND}(Cr)=(\pi+1)^{N(Cr)}<\pi^{\ome}$ of $S$ relative to a fixed sequence $Cr$ of critical formulas as in Definition \ref{df:ordind}.

In Definition \ref{df:Lp}  $M(p,n)$ is defined by $\pi^{\ome}$-recursion, 
i.e., $o(\vec{S})<\pi^{\ome}$,
and $p(Cr)$ is defined by $\veps_{\pi+1}$-recursion, i.e., $o(\vec{S})<\veps_{\pi+1}$,
which are far from $d_{\Ome}\veps_{\pi+1}<\Ome<\pi^{\ome}<\veps_{\pi+1}$.


  

 

For the moment, suppose that $M(p,n)$ 
has been defined for each $p<\ome$.
Let $\ell_{p}=M(p+1,0;Cr)$. 

If $\exi p\leq n(Cr)[\ell_{p}=\ell_{p+1}]$, then there is nothing to prove, i.e., $p(Cr)\leq n(Cr)$.
In what follows assume $\fal p\leq n(Cr)[\ell_{p}<\ell_{p+1}]$, and let $p> n(Cr)$.
Suppose $\ell_{p}<\ell_{p+1}$. Then $S^{\ell_{p}}$ is nonsolving.

Consider the number
\beqn\label{eq:mp}
m_{p}:=\max\{n<\ell_{p}: r_{n}<\pi\}
\eeqn
Note that $m_{p}$ is in the set ${\sf nd}(\vec{S}^{0,\ell_{p}})$.

\bprp\label{clm:mp}
$p>n(Cr) \spand \ell_{p}<\ell_{p+1} \Rarw \pi>r_{m_{p}}>r_{\ell_{p}}$.
\eprp
{\bf Proof} of Proposition \ref{clm:mp}.
Suppose $p>n(Cr)$, $ \ell_{p}<\ell_{p+1}$ and $r_{m_{p}}\leq r_{\ell_{p}}$.
By Lemma  \ref{lem:Lp}.\ref{lem:Lp.15} we have $M(p,n)\leq M(p+1,n)\leq M(p+1,0)=\ell_{p}$
for any $n<\ell_{p}$.
If $M(p,m_{p})=\ell_{p}$, then we would have $M(p+1,n)>\ell_{p}$ by $r_{m_{p}}\leq r_{\ell_{p}}$.
Hence $n_{p}:=M(p,m_{p})<\ell_{p}$.
On the other hand we have $r_{n_{p}}>\pi>r_{m_{p}}$, and hence 
$M(p,m_{p})< M(p+1,m_{p})\leq\ell_{p}$.
From this we see that 
$\vec{S}^{m_{p}, M(p+1,m_{p})}$ is a proper $(p+1)$-series.
Therefore a component in the decomposition of $\vec{S}^{m_{p}, M(p+1,m_{p})}$
is a proper $p$-series.
On the other side any component except the first one is an improper $p$-series
by $r_{n}>\pi\, (m_{p}<n<\ell_{p})$, Lemma \ref{lem:pseries}.\ref{lem:pseries.42} and $p\geq n(Cr)$.
Hence the first component $\vec{S}^{m_{p}, m'}\, (m'\leq M(p,m_{p}))$ is a proper $p$-series.
But then $\#\{rk(S^{n})>\pi: m_{p}<n<m'\}\geq p-1\geq n(Cr)$ by Lemma \ref{lem:pseries}.\ref{lem:pseries.42}.
This is a contradiction.
We have shown Proposition  \ref{clm:mp}.
\eprf

For each $p>n(Cr)$,
let $\vec{S}_{p+1}=\{S^{n}\}_{n<\ell_{p}}$ be the $(p+1)$-section.
Decompose $\vec{S}_{p+1}$ into 
$\vec{S}_{p+1}=\vec{S}_{p}^{1}*\cdots*\vec{S}_{p}^{l_{p}-1}*\vec{S}_{p}^{l_{p}}$ 
where the last substring $\vec{S}_{p}^{l_{p}}$ is defined to be $\vec{S}^{m_{p}, \ell_{p}}$
for the number $m_{p}$ in (\ref{eq:mp}), and
each $\vec{S}_{p}^{j}=\vec{S}^{k_{p}^{j-1}, k_{p}^{j}}\, (1\leq j<l_{p})$
for $\{(0=)k_{p}^{0},k_{p}^{1},\ldots,k_{p}^{l_{p}-1}(=m_{p})\}_{<}\incl {\sf nd}(\vec{S}_{p+1})$.
Put $k_{p}^{l_{p}}=\ell_{p}$.

Let $o(\vec{S})=o(\vec{S};0)<d_{\Ome}\ome_{n(Cr)}(\pi)$ 
denote the ordinal associated to $p$-series $\vec{S}$ 
as in  \cite{esubid}, \cite{esubMahlo}, \cite{esubidea} and \cite{esubpi2}.
Note the following fact:
\beqnarr
&& \alp<\ome_{n(Cr)}(\pi) \mbox{ for any subdiagram } d_{\sig}^{q}\alp 
\mbox{ occurring in ranks, indices and } 
\nonumber \\
&& o(\vec{S}) 
\mbox{, which appear in the H-process for } Cr
\label{eq:bndfact}
\eeqnarr

Then let $\alp_{p}^{j}=o(\vec{S}_{p}^{j})$ for $0<j\leq l_{p}-1$,
and $\alp_{p}^{0}:=d_{\Ome}\ome_{n(Cr)}(\pi)$.

We have, cf. Theorem 11.7 in \cite{esubid}
\beqn
\label{eq:alpdecimp}
\Ome>\alp_{p}^{j}>\alp_{p}^{j+1}
\eeqn

Let $\gam_{p}^{j}:=r_{k_{p}^{j+1}}\, (0\leq j< l_{p})$ and
$\Del:=\pi>\gam_{p}^{j}$.
Finally let
\[
\bet_{p} := 
\ome_{n(Cr)}(\pi) + \sum_{j<l_{p}}\Del^{\alp_{p}^{j}}\gam_{p}^{j}
.\]

\blem\label{lem:pbndpi}
$d_{\Ome}\bet_{p}\in O(\mbox{{\rm T}})$,
 and
$d_{\Ome}\bet_{p}>d_{\Ome}\bet_{p+1}$ 
if $\ell_{p}<\ell_{p+1}$ and $p>n(Cr)$.
\elem
\bprf
To show $d_{\Ome}\bet_{p}\in O(\mbox{{\rm T}})$, we have to verify a condition
$
\calb_{>\Ome}(\bet_{p})<\bet_{p}
$
for a set $\calb_{>\Ome}(\bet_{p})$ of subdiagrams of $\bet_{p}$.
This is seen from (\ref{eq:bndfact}) and $\calb_{>\Ome}(\bet_{p})=\calb_{>\Ome}(\{\gam_{p}^{j}: j<l_{p}\})$.

Assume $\ell_{p}<\ell_{p+1}$ and $p>n(Cr)$.
Let
\[
J=\min\{j\leq l_{p}-2 : \gam_{p}^{j}>\gam_{p}^{l_{p}-1}\}
.\]
Note that by Proposition  \ref{clm:mp}  we have $r_{m_{p}}=\gam_{p}^{l_{p}-2}>\gam_{p}^{l_{p}-1}=r_{\ell_{p}}$.

Then as in the proof of Lemma \ref{lem:pbnd} we see
$\alp_{p+1}^{j}=\alp_{p}^{j}$ for $0\leq j\leq J$, 
$\gam_{p+1}^{j}=\gam_{p}^{j}$ for $0\leq j<J$.
Moreover we have
$\gam_{p}^{J}>\gam_{p}^{l_{p}-1}=\gam_{p+1}^{J}$.
Therefore $\bet_{p}>\bet_{p+1}$.

It remains to show $K_{\Ome}\bet_{p+1}<d_{\Ome}\bet_{p}$.
By (\ref{eq:alpdecimp}) it suffices to see 
$K_{\Ome}\{\gam_{p+1}^{j}: j<l_{p+1}\}<d_{\Ome}\bet_{p}$.

If T=$\Phi$-FIX for $\Phi=[\Pi^{0}_{1},\Pi^{0}_{1}], \Pi^{0}_{2}$,
then there is nothing to prove, i.e., 
$K_{\Ome}\{\gam_{p+1}^{j}: j<l_{p+1}\}=\emptyset$.

Consider the case $\Phi=\Pi^{0}_{1}$ and $\pi=\Ome$.
Then 
we have
$K_{\Ome}\{\gam_{p+1}^{j}: j<l_{p+1}\}<\alp_{p}^{0}=d_{\Ome}\ome_{n(Cr)}(\Ome)<d_{\Ome}\bet_{p}$
by (\ref{eq:bndfact}).
\eprf

It remains to define $M(p,n)$ by $d_{\Ome}\veps_{\pi+1}$-recursion.
This is seen from the following lemma.
 
 \blem\label{lem:opseriespi}
  Let $\vec{S}^{i}=\vec{S}^{m_{i}, k_{i}}\, (i=0,1)$ be two consecutive $p$-series, $k_{0}=m_{1}$.
 Assume $r_{m_{1}}\geq r_{m_{0}}$.
Then
$o(\vec{S}^{1})<o(\vec{S}^{0})$.
\elem
\bprf
This is seen as in Lemma \ref{lem:opseries}, i.e., Theorems 10.8 and 11.7 in \cite{esubid}.
\eprf

Thus we have shown that both $M(p,n;Cr)$ and  $p(Cr)$ are defined by $d_{\Ome}\veps_{\pi+1}$-recursion.
This yields a solution $S^{H}$ of $Cr$ for $H=M(p(Cr),0;Cr)$.

\end{document}